\documentclass[12pt,a4paper]{amsart}
\usepackage[utf8]{inputenc}	
\usepackage[T1]{fontenc}
\usepackage[in,headings]{fullpage}
\usepackage{amsbsy, amscd, amsfonts, amsmath, amsrefs, amssymb, amsthm}
\usepackage{graphicx}
\usepackage{float}
\usepackage{color}   
\usepackage[colorlinks=true,linkcolor=blue,citecolor=blue,urlcolor=blue]{hyperref}
\usepackage{comment}
\excludecomment{A}

\newtheorem{theorem}{Theorem}

\newtheorem{proposition}[theorem]{Proposition}
\newtheorem{lemma}[theorem]{Lemma}

\theoremstyle{definition}
\newtheorem{definition}[theorem]{Definition}
\newtheorem{remark}[theorem]{Remark}

\theoremstyle{remark}

\newcommand{\C}{\mathbf{C}}

\newcommand{\R}{\mathbf{R}}

\renewcommand{\Re}{\mathop{\mathrm{Re}}\nolimits}
\renewcommand{\Im}{\mathop{\mathrm{Im}}\nolimits}

\newfont{\cmbsy}{cmbsy10}
\newfont{\cmmib}{cmmib10}
\newcommand{\Orden}{\mathop{\hbox{\cmbsy O}}\nolimits}

\def\oeis#1{\href{http://oeis.org/#1}{#1}}

\begin{document}

\title{Integral Representation for Riemann-Siegel $Z(t)$ function.}
\author[Arias de Reyna]{J. Arias de Reyna}
\address{%
Universidad de Sevilla \\ 
Facultad de Matem\'aticas \\ 
c/Tarfia, sn \\ 
41012-Sevilla \\ 
Spain.} 




\email{arias@us.es, ariasdereyna1947@gmail.com}


\begin{abstract}
We apply Poisson formula for a strip to give a representation of $Z(t)$ by means of an integral. \[F(t)=\int_{-\infty}^\infty \frac{h(x)\zeta(4+ix)}{7\cosh\pi\frac{x-t}{7}}\,dx, \qquad Z(t)=\frac{\Re F(t)}{(\frac14+t^2)^{\frac12}(\frac{25}{4}+t^2)^{\frac12}}.\] After that we get the estimate \[Z(t)=\Bigl(\frac{t}{2\pi}\Bigr)^{\frac74}\Re\bigl\{e^{i\vartheta(t)}H(t)\bigr\}+O(t^{-3/4}),\] with \[H(t)=\int_{-\infty}^\infty\Bigl(\frac{t}{2\pi}\Bigr)^{ix/2}\frac{\zeta(4+it+ix)}{7\cosh(\pi x/7)}\,dx=\Bigl(\frac{t}{2\pi}\Bigr)^{-\frac74}\sum_{n=1}^\infty \frac{1}{n^{\frac12+it}}\frac{2}{1+(\frac{t}{2\pi n^2})^{-7/2}}.\] We explain how the study of this function can lead to information about the zeros of the zeta function on the critical line.
\end{abstract}

\maketitle

\section{Introduction}

Riemann stated that the number $N_0(T)$ of zeros of $\zeta(s)$ on the critical line with ordinates between $0$ and $T$ is approximately equal to the total number of zeros $N(T)$ in the critical strip. In a letter to Weierstrass he said that this assertion follows from a new expansion of the function $\Xi(t)$ that he has not sufficiently simplified to communicate. Nobody knows what expansion  Riemann is talking about. Siegel considered the possibility that it was his integral expansion studied in \cite{Si} but concluded that probably Riemann was wrong about the possibilities of this expansion.  This work is an attempt to obtain a development following ideas available to Riemann. 

It is somewhat surprising that the approximation obtained here is very similar to the one studied in our paper \cite{A60}, obtained by an entirely different method.

\section{Integral representation for \texorpdfstring{$Z(t)$}{Z(t)}.}

In this section, we define a real analytic function $F\colon\R\to\C$ such that 
its real part is equal to $Z(t)\sqrt{\frac14+t^2}\sqrt{\frac{25}{4}+t^2}$. For this we will apply the solution of the Dirichlet problem for a strip.

The solution of Dirichlet problem in a strip was considered in \cite{HIP} 
we quote its solution from Rademacher \cite{MR0364103}*{\S\ 31}.
\begin{theorem}\label{T:Dir}
The solution $u(\sigma,t)$  of the Dirichlet problem on the strip 
$S(a,b)=\{s\in\C\colon a\le \sigma\le b\}$ with boundary values 
$u(a,t)=A(t)$ and $u(b,t)=B(t)$, where $A(t)$ and $B(t)=\Orden(e^{k|t|})$ with 
$0<k<\pi/(b-a)$  is given by 
\begin{multline}
U(\sigma,t)=\frac{1}{2(b-a)}\int_{-\infty}^{+\infty}A(x)
\omega\Bigl(\frac{\sigma-a}{b-a},\frac{x-t}{b-a}\Bigr)\,dx\\+
\frac{1}{2(b-a)}\int_{-\infty}^{+\infty}B(x)
\omega\Bigl(\frac{b-\sigma}{b-a},\frac{x-t}{b-a}\Bigr)\,dx.
\end{multline}
where 
\begin{equation}
\omega(\sigma,t)=\frac{\sin\pi\sigma}{\cosh\pi t-\cos\pi\sigma}.
\end{equation}
\end{theorem}

\begin{theorem}
There is a holomorphic function $f$  in the strip $S(-3,5)$ such that for $\frac12<\sigma<5$ and $t\in\R$ 
\begin{equation}\label{E:new1}
Z(t)\sqrt{\tfrac14+t^2}\sqrt{\tfrac{25}{4}+t^2}=\Re\Bigl\{\frac{1}{2\sigma-1}\int_{-\infty}^\infty 
\frac{f(\sigma+ix)}{\cosh\pi\frac{x-t}{2\sigma-1}}\,dx\Bigr\}.
\end{equation}
\end{theorem}
\begin{proof}
Consider the function $\Phi(s)=(s+2)s(1-s)(3-s)\zeta(s)\zeta(1-s)$. It is clear that 
$\Phi(s)$ is an entire function that satisfies the functional equation $\Phi(s)=\Phi(1-s)$. The function $\Phi(s)$ have simple zeros at $-4$, $-6$ \dots\ and at  $5$, $7$ \dots\  due to the trivial zeros of $\zeta(s)$. The zeros on the strip $S(-4,5)$ are all of even order, at $-2$, $3$ and at each nontrivial zero $\rho$ of $\zeta(s)$. 

In the critical line, we have $\zeta(\frac12+it)=Z(t)e^{-i\vartheta(t)}$ where $Z(t)$ and $\vartheta(t)$ real functions that  are real analytic, and such that 
\begin{equation}
Z(t)=Z(-t), \qquad \vartheta(-t)=-\vartheta(t).
\end{equation}
 It follows that 
\begin{equation}
\Phi(\tfrac12+it)=(\tfrac14+t^2)(\tfrac{25}{4}+t^2)Z(t)^2.
\end{equation}

Since the strip $S(-3,5)$ is simply connected and the zeros of $\Phi(s)$ are of even order,  there is an analytic function $f\colon S(-4,5)\to\C$ with $f(s)^2=\Phi(s)$. We fix  the sign of $f(s)$ so that 
\begin{equation}\label{E:fcritic}
f(\tfrac12+it)=Z(t)\sqrt{\tfrac14+t^2}\sqrt{\tfrac{25}{4}+t^2},\qquad \text{for $t\in\R$},
\end{equation}
where we take $\sqrt{\tfrac14+t^2}>0$ and $\sqrt{\tfrac{25}{4}+t^2}>0$.   
Notice that $Z(0)=\zeta(1/2)=-1.46 <0$, so that $f(\frac12)<0$.

We apply Theorem \ref{T:Dir} to the harmonic function $f(s)=f(\sigma+it)$ on the strip $S(1-\sigma,\sigma)$. The boundary values are $A(t)=f(1-\sigma+it)$ and $B(t)=f(\sigma+it)$. Since $f(s)$ is real on the real axis and on the critical line, we have $f(1-\sigma+it)=\overline{f(\sigma+it)}$.  Therefore, we obtain 
\begin{align*}
f(\tfrac12+it)&=Z(t)\sqrt{\tfrac14+t^2}\sqrt{\tfrac{25}{4}+t^2}\\
&=\frac{1}{2(2\sigma-1)}
\int_{-\infty}^{\infty} A(x)\omega\bigl(\tfrac{\sigma-\frac12}{2\sigma-1},\tfrac{x-t}{2\sigma-1}\bigr)\,dx+
\frac{1}{2(2\sigma-1)}\int_{-\infty}^{+\infty}B(x)
\omega\bigl(\tfrac{\sigma-\frac12}{2\sigma-1},\tfrac{x-t}{2\sigma-1}\bigr)\,dx.
\end{align*}
The two integrals are complex conjugate numbers; therefore,
\[Z(t)\sqrt{\tfrac14+t^2}\sqrt{\tfrac{25}{4}+t^2}=2\Re\frac{1}{2(2\sigma-1)}
\int_{-\infty}^{\infty} B(x)\omega\bigl(\tfrac12,\tfrac{x-t}{2\sigma-1}\bigr)\,dx.\]
We have 
\[\omega\bigl(\tfrac12,\tfrac{x-t}{2\sigma-1}\bigr)=\frac{\sin\frac{\pi}{2}}
{\cosh\pi\frac{x-t}{2\sigma-1}-\cos\frac{\pi}{2}}=
\frac{1}{\cosh\pi\frac{x-t}{2\sigma-1}},\]
and the result follows.
\end{proof}

There are several interesting cases of \eqref{E:new1}, but we start by trying to 
obtain the consequences of case $\sigma=4$. 
\[Z(t)\sqrt{\tfrac14+t^2}\sqrt{\tfrac{25}{4}+t^2}=\Re\Bigl\{\frac{1}{7}\int_{-\infty}^\infty 
\frac{f(4+ix)}{\cosh\pi\frac{x-t}{7}}\,dx\Bigr\}.\]
Therefore, we define for $t\in \R$
\begin{equation}\label{E:intform1}
F(t)=\int_{-\infty}^\infty 
\frac{f(4+ix)}{7\cosh\pi\frac{x-t}{7}}\,dx.
\end{equation}
It is easy to prove that $F(t)$ is analytic for $|\Im(t)|<\frac72$.

\subsection{The function  \texorpdfstring{$f(4+ix)$}{f(4+ix)}.}

\begin{proposition}
The function $f(4+ix)=h(x)\zeta(4+ix)$ where $h(x)$ is the real analytic function 
defined by 
\begin{equation}\label{E:defh}
h(x)=\frac{\sqrt{2}}{(2\pi)^2}
\Bigl\{(6+ix)(4+ix)(3+ix)(1+ix)(2\pi)^{-ix}\cosh\tfrac{\pi x}{2}\Gamma(4+ix)\Bigr\}^{1/2},\end{equation}
that satisfies $h(0)=\frac{3}{\pi^2}\sqrt{6}>0$.
\end{proposition}

\begin{proof}
For $-2<\sigma<3$ the function $\Phi(\sigma)$ is continuous and does not vanish,
so its sign coincides with the sign at $s=1/2$. Therefore $\Phi(\sigma)>0$ in 
this range. In the same range $f(\sigma)<0$ because we took $f(\frac12)<0$ (take $t=0$ in 
\eqref{E:fcritic}).

$f(s)$ is defined on the strip $S(-4,5)$ as the analytic square root of the 
function $\Phi(s)$ that is negative at $s=\frac12$.

By the functional equation of $\zeta(s)$ we have 
\begin{align*}
\Phi(s)&=(s+2)s(1-s)(3-s)\zeta(s)\zeta(1-s)\\
&=2(s+2)s(1-s)(3-s)(2\pi)^{-s}\cos\tfrac{\pi s}{2}
\Gamma(s) \zeta(s)^2.
\end{align*}
Notice that the function $2(s+2)s(1-s)(3-s)(2\pi)^{-s}\cos\tfrac{\pi s}{2}\Gamma(s)$ is analytic in $s=4+ix$ with $x\in\R$ and does not vanish so that there is a real analytic function 
$h\colon\R\to\C$  defined by \eqref{E:defh}. By a simple computation we obtain 
\[\frac{\Phi(4+ix)}{\zeta(4+ix)^2}=\frac{2}{(2\pi)^4}(6+ix)(4+ix)(3+ix)(1+ix)
(2\pi)^{-ix}\cosh\tfrac{\pi x}{2}\Gamma(4+ix).\]
From this we see that $h(0)=\pm \frac{3}{\pi^2}\sqrt{6}$. We may define $h(x)$ as the only continuous square root with $h(0)=\frac{3}{\pi^2}\sqrt{6}$.

With these elections we have $f(4+ix)=\pm h(x)\zeta(4+ix)$.  We know that 
$f(x)<0$ for $-2<x<3$, and $f(4)>0$, because $f$ changes sign at $x=3$. 
Therefore, taking $h(0)>0$ we obtain the equality 
$f(4+ix)=h(x)\zeta(4+ix)$. 
\end{proof}

With the above definitions, we may write \eqref{E:intform1} as 
\begin{equation}\label{E:intform2}
F(t)=\int_{-\infty}^\infty 
\frac{h(x)\zeta(4+ix)}{7\cosh\pi\frac{x-t}{7}}\,dx, \qquad Z(t)=\frac{\Re F(t)}{(\frac14+t^2)^{\frac12}(\frac{25}{4}+t^2)^{\frac12}}.
\end{equation}

\begin{proposition}
There are real analytic real functions $\rho\colon\R\to\R$ and $\alpha\colon\R\to\R$ such that for all $x\in\R$ we have $h(x)=\rho(x)e^{i\alpha(x)}$.
\end{proposition}

\begin{proof}
Since $h(x)\ne0$, we may use the result in \cite{AL}*{Cor.2.3} (see also Proposition \ref{P:1} in Section \ref{appendix}) to show that there are a real analytic function $\alpha(x)$ such that $h(x)=|h(x)|e^{i\alpha(x)}$. The function $\rho(x):=|h(x)|$ is real analytic, because $h(x)$ does not vanish.

Since $\rho(0)>0$  and $h(0)>0$ we have $\alpha(0)=2k\pi$. For definiteness, we take $\alpha(0)=0$. 
\end{proof}

\begin{proposition}
For $x\to+\infty$, we have the two asymptotic expansions
\begin{equation}
\begin{aligned}
\log\rho(x)\sim2\log x+\frac74\log\frac{x}{2\pi}-\sum_{n=1}^\infty\frac{(-1)^n(6^{2n}+4^{2n}+2\cdot 3^{2n}+2^{2n}+2)}{4n x^{2n}},\\
\end{aligned}
\end{equation}
\begin{equation}
\begin{aligned}
\alpha(x)&\sim\frac{x}{2}\log\frac{x}{2\pi}-\frac{x}{2}+\frac{15\pi }{8}\\&-
\sum_{n=0}^\infty \frac{(-1)^n(3\cdot 6^{2n}+2\cdot 4^{2n}+3\cdot 3^{2n}+2^{2n}+1)}{(2n+1) x^{2n+1}}+\sum_{n=1}^\infty 
\frac{(-1)^nB_{2n}}{4n(2n-1)x^{2n-1}}.
\end{aligned}
\end{equation}
\end{proposition}

\begin{proof}
Consider the function between parenthesis in \eqref{E:defh}
\[(6+ix)(4+ix)(3+ix)(1+ix)(2\pi)^{-ix}\cosh\tfrac{\pi x}{2}\Gamma(4+ix),\]
it is equal to 
\[(6+ix)(4+ix)(3+ix)^2(2+ix)(1+ix)^2(2\pi)^{-ix}\cosh\tfrac{\pi x}{2}\Gamma(1+ix).\]
This function is equal to $2^43^3$ for $x=0$ and its factors are all positive at $x=0$, so we consider its logarithm as the sum of logarithms that are real at $x=0$.
\begin{multline*}
\log(6+ix)+\log(4+ix)+2\log(3+ix)+\log(2+ix)+2\log(1+ix)\\+\log\bigl(\cosh\tfrac{\pi x}{2}\bigr)-ix\log(2\pi)+\log\Gamma(1+ix).
\end{multline*}
Each term of the above sum has a well-known asymptotic expansion for $x\to+\infty$
\begin{gather*}
\log(6+ix)=\log x+\frac{\pi i}{2}-\sum_{n=1}^\infty\frac{1}{n}\Bigl(\frac{6i}{x}\Bigr)^n,\quad
\log(4+ix)=\log x+\frac{\pi i}{2}-\sum_{n=1}^\infty\frac{1}{n}\Bigl(\frac{4i}{x}\Bigr)^n,\\
2\log(3+ix)=2\log x+\pi i-\sum_{n=1}^\infty\frac{2}{n}\Bigl(\frac{3i}{x}\Bigr)^n,\quad
\log(2+ix)=\log x+\frac{\pi i}{2}-\sum_{n=1}^\infty\frac{1}{n}\Bigl(\frac{2i}{x}\Bigr)^n,\\
2\log(1+ix)=2\log x+\pi i-\sum_{n=1}^\infty\frac{2}{n}\Bigl(\frac{i}{x}\Bigr)^n,\\
\log\bigl(\cosh\tfrac{\pi x}{2}\bigr)=-\log 2+\frac{\pi x}{2}+\log(1+e^{-\pi x})
\sim \frac{\pi x}{2}-\log 2,
\end{gather*}
\begin{multline*}
-ix\log(2\pi)+\log\Gamma(1+ix)\\ \sim-\frac{\pi x}{2}+\frac{1}{2}\log\frac{x}{2\pi}+
\log(2\pi)
+ix\log\frac{x}{2\pi}-ix+\frac{\pi i}{4}+i\sum_{n=1}^\infty 
\frac{(-1)^nB_{2n}}{2n(2n-1)x^{2n-1}}.
\end{multline*}
Then $\log h(x)$ is easily obtained from these values. From \eqref{E:defh}
we have to add $\frac12$ of the above expansions and add the log of $\sqrt{2}/(2\pi)^2$.
According to our election $\alpha(0)=0$ we must take the logarithm 
of $\sqrt{2}/(2\pi)^2$ equal to $\log(\sqrt{2}/(2\pi)^2)$ the real logarithm. Therefore,
\begin{align*}
\log h(x)&=\frac12\log 2-2\log(2\pi)\\
&+\frac72\log x+\frac{7\pi i}{4}-\sum_{n=1}^\infty
\frac{6^n+4^n+2\cdot 3^n+2^n+2}{2n}\Bigl(\frac{i}{x}\Bigr)^n\\
&+\frac{\pi x}{4}-\frac12\log 2-\frac{\pi x}{4}+\frac14\log\frac{x}{2\pi}+\frac12\log(2\pi)\\
&+\frac{ix}{2}\log\frac{x}{2\pi}-\frac{ix}{2}+\frac{\pi i}{8}+i\sum_{n=1}^\infty 
\frac{(-1)^nB_{2n}}{4n(2n-1)x^{2n-1}}.
\end{align*}
Separating the real and imaginary parts, we 
obtain the asymptotic expansion of $\log\rho(x)$ and $\alpha(x)$ 
\begin{multline*}
\alpha(x)\sim\frac{x}{2}\log\frac{x}{2\pi}-\frac{x}{2}+\frac{15\pi }{8}\\-
\sum_{n=0}^\infty \frac{(-1)^n(3\cdot 6^{2n}+2\cdot 4^{2n}+3\cdot 3^{2n}+2^{2n}+1)}{(2n+1) x^{2n+1}}+\sum_{n=1}^\infty 
\frac{(-1)^nB_{2n}}{4n(2n-1)x^{2n-1}}.
\end{multline*}
\[\log\rho(x)\sim -\frac32\log(2\pi)+\frac72\log x+\frac14\log \frac{x}{2\pi}-
\sum_{n=1}^\infty\frac{(-1)^n(6^{2n}+4^{2n}+2\cdot 3^{2n}+2^{2n}+2)}{4n x^{2n}}.\]
or by summing the two convergent series we obtain 
\[\log\rho(x)\sim2\log x+\frac74\log\frac{x}{2\pi}+\frac14\log\frac{(1+x^2)^2(4+x^2)(16+x^2)(9+x^2)^2(36+x^2)}{x^{14}}
.\]
The difference between these two functions is exponentially small.

\begin{multline*}
\alpha(x)\sim\frac{x}{2}\log\frac{x}{2\pi}-\frac{x}{2}+\frac{15\pi }{8}-
\arctan\frac{1}{x}-\frac12\arctan\frac{2}{x}-\arctan\frac{3}{x}\\
-\frac12\arctan\frac{4}{x}-\frac12\arctan\frac{6}{x}
+\sum_{n=1}^\infty 
\frac{(-1)^nB_{2n}}{4n(2n-1)x^{2n-1}}.
\end{multline*}
\end{proof}

\begin{lemma}\label{lrhofirstterms}
The first terms of the expansion of $\log\rho(x)$ are given by
\[\log\rho(x)=-\frac74\log(2\pi)+\frac{15}{4}\log x+\frac{19}{x^2}-\frac{433}{2x^4}+\frac{13069}{3x^6}-\frac{439633}{4x^8}+\cdots\]
\end{lemma}
Therefore, we obtain the next Lemma.
\begin{lemma}\label{rhobound}
We have 
\[\rho(x)=x^2\Bigl(\frac{x}{2\pi}\Bigr)^{7/4}\Bigl(1+\frac{19}{x^2}+\Orden(x^{-4})\Bigr),\qquad x\to+\infty.\]
\end{lemma}
Analogously, the first terms of the asymptotic expansion of $\alpha(x)$ by the above reasoning are equal to 
\[\alpha(x)=\frac{x}{2}\log\frac{x}{2\pi}-\frac{x}{2}+\frac{15\pi}{8}-\frac{241}{24x}+\frac{41279}{720x^3}-\frac{2348641}{2520x^5}+\cdots\]
Therefore,
\begin{lemma}\label{alphaexp}
We have for $x\to+\infty$
\[\alpha(x)=\frac{x}{2}\log\frac{x}{2\pi}-\frac{x}{2}+\frac{15\pi}{8}-\frac{241}{24x}
+\Orden(x^{-3}).\]
\end{lemma}

\section{First approximation to \texorpdfstring{$Z(t)$}{Z(t)}.}

Let us define  
\begin{equation}\label{deftheta}
\theta(x)=\frac{x}{2}\log\frac{x}{2\pi}-\frac{x}{2}+\frac{15\pi}{8}-\frac{241}{24x},
\qquad x>0,
\end{equation}
\begin{equation}\label{defrho0}
\rho_0(x)=x^2\Bigl(\frac{x}{2\pi}\Bigr)^{7/4}\Bigl(1+\frac{19}{x^2}\Bigr)=
(2\pi)^{-7/4}x^{7/4}(19+x^2),\qquad x>0,
\end{equation}
and 
\[L_1(x,t)=1+\frac{15x}{4t}+\frac{ix^2}{4t}+\frac{165 x^2}{32t^2}
+\frac{241ix}{24t^2}+\frac{41ix^3}{48t^2}-\frac{x^4}{32t^2}.\]

The object of this section is to prove the following theorem
\begin{theorem}
For $t\in\R$ let
\begin{equation}
G(t)=e^{i\theta(t)}\rho_0(t)\int_{-\infty}^\infty L_1(x,t)\Bigl(\frac{t}{2\pi}\Bigr)^{ix/2}
\frac{\zeta(4+it+ix)}{7\cosh\frac{\pi x}{7}}\,dx.
\end{equation}
Then $G(t)$ is real analytic and for $t\to+\infty$ we have $F(t)=G(t)+\Orden(t^{3/4}\log^6t)$, and consequently
\begin{equation}\label{E:1staprox}
Z(t)=\frac{\Re(G(t))}{\sqrt{\frac14+t^2}\sqrt{\frac{25}{4}+t^2}}+\Orden(t^{-5/4}\log^6 t).
\end{equation}
\end{theorem}
\begin{proof}
In the next Lemmas and Propositions we 
define the functions $F_1(t)$,  $F_2(t)$, $F_3(t)$ and $F_4(t)$, such that 
$F_4(t)=G(t)$ and $F(t)=F_1(t)+\Orden(t^{-1/4})$, $F_1(t)=F_2(t)+\Orden(t^{3/4})$,
$F_2(t)=F_3(t)+\Orden(t^{3/4}\log^{6}t)$ and $F_3(t)-F_4(t)=\Orden(t^{-1/4})$. 
And the Theorem follows from this and \eqref{E:intform2}.
\end{proof}

\begin{lemma}[Gabcke]\label{Gabcke}
The incomplete Gamma functions satisfies the inequality
\[\Gamma(a,x):=\int_x^\infty v^{a-1}e^{-v}\,dv\le ae^{-x}x^{a-1},\qquad x>a\ge1.\]
\end{lemma}
\begin{proof}
See \cite{G}*{p.145, p.84 in the electronic version}.
\end{proof}

\begin{lemma}\label{L:instaurate}
Let $\alpha=28/\pi$ and
\[F_1(t)=\int_{t-\alpha\log t}^{t+\alpha\log t}\frac{h(x)\zeta(4+ix)}{7\cosh\frac{\pi}{7}(x-t)}\,dx.\]
Then for $t\to+\infty$ we have $F(t)=F_1(t)+\Orden(t^{-1/4})$. \end{lemma}
\begin{proof}
The  difference $F(t)-F_1(t)$ is the sum of two integrals, we bound each one.
Let $b>0$, we shall take at the end $b=\alpha\log t$. By Lemma \ref{rhobound} $\rho(x)\le C x^{15/4}$ for $x\ge x_0$. Also 
$|\zeta(4+ix)|$ is bounded for $x$ real. Therefore, we have for $t\to+\infty$
\begin{align*}
\Bigl|\int_{t+b}^\infty\frac{h(x)\zeta(4+ix)}{7\cosh\frac{\pi}{7}(x-t)}\,dx\Bigr|
&\ll 
\int_{t+b}^\infty\frac{x^{15/4}}{7\cosh\frac{\pi}{7}(x-t)}\,dx\\
&=
\int_{b/t}^\infty\frac{t^{15/4}(1+y)^{15/4}}{7\cosh\pi y t/7}t\,dy.
\end{align*}
We divide the integral into two other (assuming $t>b$)
\begin{align*}
&=t^{15/4}\int_{b/t}^1\frac{(1+y)^{15/4} t}{7\cosh\pi y t/7}\,dy
+t^{15/4}\int_{1}^\infty\frac{(1+y)^{15/4}t}{7\cosh\pi y t/7}\,dy\\
&\le2^{15/4}t^{15/4}\int_{b/t}^1\frac{t}{\cosh\pi y t/7}\,dy
+2^{15/4}t^{15/4}\int_{1}^\infty\frac{y^{15/4}t}{\cosh\frac{\pi y t}{7}}\,dy\\
&\ll t^{15/4}\int_{b/t}^\infty t e^{-\frac{\pi y t}{7}}\,dy
+t^{15/4}\int_{1}^\infty y^{15/4} t e^{-\frac{\pi y t}{7}}\,dy\\
&=\frac{7}{\pi}t^{15/4} e^{-\pi b/7}+t^{15/4}\int_{\pi t/7}^\infty\Bigl(\frac{7x}{\pi t}\Bigr)^{15/4}te^{-x}\frac{7}{\pi t}\,dx\\
&\ll t^{15/4} e^{-\pi b/7}+\int_{\pi t/7}^\infty x^{15/4}e^{-x}\,dx
\end{align*}
The last integral can be bounded by Gabcke's lemma. In the first one we take as 
$b=\alpha\log t$ and  obtain 
\[\ll t^{15/4} t^{-4}+t^{15/4}e^{-\pi t/7}=\Orden(t^{-1/4}).\]

Now we turn to the second integral, for some positive $a<t/2$, (we will take at the end $a=\alpha\log t$)
\begin{align*}
\Bigl|\int_{-\infty}^{t-a}\frac{h(x)\zeta(4+ix)}{7\cosh\frac{\pi}{7}(x-t)}&\,dx\Bigr|
\\&\le \Bigl|\int_{t/2}^{t-a}\frac{h(x)\zeta(4+ix)}{7\cosh\frac{\pi}{7}(x-t)}\,dx\Bigr|+
\Bigl|\int_{-\infty}^{t/2}\frac{h(x)\zeta(4+ix)}{7\cosh\frac{\pi}{7}(x-t)}\,dx\Bigr|:= I_1+I_2.
\end{align*}
For the second integral here, we notice that $h(x)$ is bounded for $|x|<1$ so that we have
\begin{multline*}
I_2=\Bigl|\int_{-\infty}^{t/2}\cdots\Bigr|\ll \int_{-1}^1
\frac{1}{\cosh\frac{\pi}{7}(x-t)}\,dx\\
+\int_{1}^{t/2}\frac{x^{15/4}}{\cosh\frac{\pi}{7}(x-t)}\,dx+\int_{-\infty}^{-1}\frac{|x|^{15/4}}{\cosh\frac{\pi}{7}(x-t)}\,dx:=J_1+J_2+J_3.
\end{multline*}
It is easy to see that $J_1\ll e^{-\pi t/7}$. For $J_2$ we can obtain a better 
bound, but it is sufficient to bound the integral by the maximum of the function 
by the length of the interval 
\[J_2\ll t^{15/4}e^{-\pi t/14}\cdot  t=\Orden(t^{-1/4}).\]
For $J_3$ we reduce it to an incomplete gamma integral 
\[
J_3=\int_1^\infty \frac{x^{15/4}}{\cosh\frac{\pi}{7}(x+t)}\,dx\ll
\int_1^\infty x^{15/4}e^{-\pi\frac{x+t}{7}}\,dx\ll e^{-\pi t/7}.
\]
All is reduced to the bound of $I_1$
\begin{align*}
I_1&=\Bigl|\int_{t/2}^{t-a}\frac{h(x)\zeta(4+ix)}{7\cosh\frac{\pi}{7}(x-t)}\,dx\Bigr|
\ll \int_{t/2}^{t-a}\frac{x^{15/4}}{\cosh\frac{\pi}{7}(x-t)}\,dx\\
&\ll t^{15/4}\int_{t/2}^{t-a}\frac{dx}{\cosh\frac{\pi}{7}(x-t)}\le 
t^{15/4}\int_{-\infty}^{-a}\frac{dy}{\cosh\frac{\pi y}{7}}\\
&\le 2 t^{15/4}\int_{-\infty}^{-a}e^{\frac{\pi y}{7}}\,dy\ll 
t^{15/4} e^{-\pi a/7}.
\end{align*}
We see that this is $\Orden(t^{-1/4})$ if we choose $a=\alpha\log t$ with $\alpha=28/\pi$,
as we wanted.
\end{proof}

\begin{lemma}
For $t-\alpha\log t>0$ let us define
\[F_2(t)=\int_{t-\alpha\log t}^{t+\alpha\log t}\frac{\rho_0(x)e^{i\theta(x)}\zeta(4+ix)}{7\cosh\frac{\pi}{7}(x-t)}\,dx,\]
where $\rho_0(x)$ and $\theta(x)$ are defined in \eqref{defrho0} and \eqref{deftheta}.
Then $F_1(t)=F_2(t)+\Orden(t^{3/4})$.
\end{lemma}
\begin{proof} By Lemma \ref{rhobound} and the definition of $\rho_0(x)$ we have for $x\to+\infty$
\[\rho(x)=\rho_0(x)+\Orden(x^{-1/4}),\qquad \rho_0(x)=\Orden(x^{15/4}).\]
By Lemma \ref{alphaexp} and the definition of $\theta(x)$ we have
\[\alpha(x)=\theta(x)+\Orden(x^{-3}).\]
Then 
\begin{align*}
h(x)&=\rho(x)e^{i\alpha(x)}\\
&=\bigl(\rho_0(x)+\Orden(x^{-1/4})\bigr) e^{i\theta(x)+\Orden(x^{-3})}\\
&=\rho_0(x)e^{i\theta(x)}\bigl(1+\Orden(x^{-4})\bigr)\bigl(1+\Orden(x^{-3})\bigr)\\
&=\rho_0(x)e^{i\theta(x)}\bigl(1+\Orden(x^{-3})\bigr)\\
&=\rho_0(x)e^{i\theta(x)}+\Orden(x^{3/4}).\\
\end{align*}
It follows that the difference $F_1(t)-F_2(t)$ is bounded by 
\[F_1(t)-F_2(t)\ll \int_{t-\alpha\log t}^{t+\alpha\log t}\frac{x^{3/4}}{7\cosh\frac{\pi}{7}(x-t)}\,dx\]
In the integration interval $x^{3/4}\ll t^{3/4}$ and we have
\[F_1(t)-F_2(t)\ll t^{3/4}\int_{t-\alpha\log t}^{t+\alpha\log t}\frac{1}{7\cosh\frac{\pi}{7}(x-t)}\,dx\ll t^{3/4},\]
since the integral on $\R$ is bounded by a constant.
\end{proof}
 
\begin{lemma}
For $t-\alpha\log t>0$, define
\[F_3(t)=e^{i\theta(t)}\rho_0(t)\int_{-\alpha\log t}^{\alpha\log t}
L_1(x,t)\Bigl(\frac{t}{2\pi}\Bigr)^{ix/2}
\frac{\zeta(4+it+ix)}{7\cosh(\pi x/7)}\,dx,\]
where
\begin{equation}\label{defL1}
L_1(x,t)=1+\frac{15x}{4t}+\frac{ix^2}{4t}+\frac{165x^2}{32t^2}
+\frac{241ix}{24t^2}+\frac{41ix^3}{48t^2}-\frac{x^4}{32t^2}.\end{equation}
Then we have $F_2(t)=F_3(t)+\Orden(t^{3/4}\log^{6}t)$.
\end{lemma}
\begin{proof}
First, notice that
\[
F_2(t)=\int_{t-\alpha\log t}^{t+\alpha\log t}\frac{\rho_0(x)
e^{i\theta(x)}\zeta(4+ix)}{7\cosh\frac{\pi}{7}(x-t)}\,dx=
\int_{-\alpha\log t}^{\alpha\log t}\frac{\rho_0(t+x)
e^{i\theta(t+x)}\zeta(4+it+ix)}{7\cosh(\pi x/7)}\,dx.
\]
Note also that $|x|\le\alpha\log t$ when $x$ is in the integration interval.
We apply Taylor's theorem to the two functions (see Rudin \cite{Ru}*{p.~110}). 

For $x$ in the integration interval there exists $\xi$ between $t$ and $x+t$,
therefore $\xi\in 
(t-\alpha\log t, t+\alpha\log t)$ with 
\begin{align*}
\theta(t+x)=&\Bigl(\frac{t}{2}\log\frac{t}{2\pi}-\frac{t}{2}
+\frac{15\pi}{8}-\frac{241}{24t}\Bigr)+\Bigl(\frac12\log\frac{t}{2\pi}
+\frac{241}{24t^2}\Bigr)x+\Bigl(\frac{1}{4t}-\frac{241}{24t^3}\Bigr)x^2\\
&-\Bigl(\frac{1}{12t^2}-\frac{241}{24t^4}\Bigr)x^3+\Bigl(\frac{1}{24\xi^3}-
\frac{241}{24\xi^5}\Bigr)x^4.
\end{align*}
So,
\begin{equation}\label{Taylortheta}
\theta(t+x)=\theta(t)+\Bigl(\frac12\log\frac{t}{2\pi}
+\frac{241}{24t^2}\Bigr)x+\frac{1}{4t}x^2
-\frac{1}{12t^2}x^3+\Orden(t^{-3}\log^4t).
\end{equation}

Recall that 
\[\rho_0(x)=x^2\Bigl(\frac{x}{2\pi}\Bigr)^{7/4}\Bigl(1+\frac{19}{x^2}\Bigr)=
(2\pi)^{-7/4}x^{7/4}(19+x^2).\]
For $t>0$ large enough and $x\in(-\alpha\log t, \beta\log t)$, 
we have Taylor's expansion with $\xi\in 
(t-\alpha\log t, t+\alpha\log t)$
\begin{align*}
(2\pi)^{7/4}\rho_0(t+x)&=t^{7/4}(19+t^2)+t^{7/4}\Bigl(\frac{15t}{4}+\frac{133}{4t}\Bigr)x\\
&+t^{7/4}\Bigl(\frac{165}{32}+\frac{399}{32t^2}\Bigr)x^2+
\xi^{7/4}\Bigl(\frac{385}{128\xi}-\frac{133}{128\xi^3}\Bigr)x^3.
\end{align*}
The last term is $\Orden(t^{3/4}(\log t)^3)$. Eliminating all terms 
of order less than or equal to this one, we obtain 
\[
(2\pi)^{7/4}\rho_0(t+x)=t^{7/4}(19+t^2)+
\frac{15 t^{11/4}}{4}x
+\frac{165 t^{7/4}}{32}x^2+\Orden(t^{3/4}(\log t)^3).
\]
So that $(2\pi)^{7/4}\rho_0(t+x)$ is equal to 
\begin{equation}\label{Taylorrho}
(2\pi)^{7/4}\rho_0(t)\Bigl(1+\frac{15t}{4(19+t^2)}x+
\frac{165}{32(19+t^2)}x^2+\Orden(t^{-3}(\log t)^3)\Bigr).
\end{equation}

By \eqref{Taylortheta} and \eqref{Taylorrho} we have
\begin{equation}\label{E:F2F3}
(2\pi)^{7/4}\rho_0(t+x)e^{i\theta(t+x)}=(2\pi)^{7/4}\rho_0(t)e^{i\theta(t)}
\Bigl(\frac{t}{2\pi}\Bigr)^{ix/2}K(x,t),
\end{equation}
where
\begin{multline*}
K(x,t)=\exp\Bigl\{i\Bigl(\frac{241 x}{24t^2}+\frac{x^2}{4t}
-\frac{x^3}{12t^2}+\Orden(t^{-3}\log^4t)\Bigr)\Bigr\}\cdot\\
\cdot\Bigl(1+\frac{15t}{4(19+t^2)}x+
\frac{165}{32(19+t^2)}x^2+\Orden(t^{-3}(\log t)^3)\Bigr).
\end{multline*}
The exponential factor is equal to
\[\exp\Bigl(\frac{241 i x}{24t^2}+\frac{ix^2}{4t}-\frac{i x^3}{12t^2}
\Bigr)(1+\Orden(t^{-3}\log^4t)).\]
Expanding  the exponential and retaining only  terms greater than $t^{-3}$
yields 
\begin{multline*}
\Bigl(1+\frac{ix^2}{4t}+\frac{241ix}{24t^2}-\frac{ix^3}{12t^2}-\frac{x^4}{32t^2}
+\Orden(t^{-3}(\log t)^6)\Bigr)(1+\Orden(t^{-3}\log^4t))\\
=1+\frac{ix^2}{4t}+\frac{241ix}{24t^2}-\frac{ix^3}{12t^2}-\frac{x^4}{32t^2}
+\Orden(t^{-3}(\log t)^6).
\end{multline*}
Expanding also the second factor in the expression for $K(x,t)$ and retaining
only the terms greater than $t^{-3}$, yields
\[1+\frac{15x}{4t}+\frac{165x^2}{32t^2}+\Orden(t^{-3}(\log t)^3).\]
Therefore,
\begin{align*}
K(x,t)&=\Bigl(1+\frac{ix^2}{4t}+\frac{241ix}{24t^2}-\frac{ix^3}{12t^2}-\frac{x^4}{32t^2}
+\Orden(t^{-3}(\log t)^6)\Bigr)\cdot\\
&\mskip120mu\cdot
\Bigl(1+\frac{15x}{4t}+\frac{165 x^2}{32t^2}+\Orden(t^{-3}(\log t)^3)\Bigr)\\
&=1+\frac{15x}{4t}+\frac{ix^2}{4t}+\frac{165x^2}{32t^2}
+\frac{241ix}{24t^2}+\frac{41ix^3}{48t^2}-\frac{x^4}{32t^2}+\Orden(t^{-3}(\log t)^6).
\end{align*}
That is $K(x,t)=L_1(x,t)+\Orden(t^{-3}(\log t)^6)$ and therefore  \eqref{E:F2F3}
can be written as
\[
\rho_0(t+x)e^{i\theta(t+x)}=\rho_0(t)e^{i\theta(t)}
\Bigl(\frac{t}{2\pi}\Bigr)^{ix/2}(L_1(x,t)+\Orden(t^{-3}(\log t)^6)),
\]
It follows that the difference $F_2(t)-F_3(t)$ is bounded by 
\[F_2(t)-F_3(t)\ll t^{15/4}\cdot t^{-3}(\log t)^6
\int_{-\alpha\log t}^{\alpha\log t}\frac{|\zeta(4+it+ix)|}{7\cosh(\pi x/7)}\,dx
=\Orden(t^{3/4}(\log t)^6),\]
since the integral is bounded by a constant.
\end{proof}

To simplify $L_1(x,t)$ we need to consider again the complete integral.
\begin{lemma}\label{L:restaurate}
Let 
\[G(t)=F_4(t)=e^{i\theta(t)}\rho_0(t)\int_{-\infty}^{+\infty}
L_1(x,t)\Bigl(\frac{t}{2\pi}\Bigr)^{ix/2}
\frac{\zeta(4+it+ix)}{7\cosh(\pi x/7)}\,dx,\]
Then $F_3(t)-F_4(t)=\Orden(t^{-5/4}\log t)$.
\end{lemma}
\begin{proof}
We need only to bound  the  integrals extended to the intervals 
$(-\infty,-\alpha\log t)$ and $(\alpha\log t,+\infty)$.  
Since $|\zeta(4+ix)|\le\zeta(4)$ and $\rho_0(t)\ll t^{15/4}$ we have
for the first integral
\[
I_1:=e^{i\theta(t)}\rho_0(t)\int_{-\infty}^{-\alpha\log t}
L_1(x,t)\Bigl(\frac{t}{2\pi}\Bigr)^{ix/2}
\frac{\zeta(4+it+ix)}{7\cosh(\pi x/7)}\,dx\ll t^{15/4}
\int_{-\infty}^{-\alpha\log t}
\frac{|L_1(x,t)|}{\cosh(\pi x/7)}\,dx.
\]
The definition \eqref{defL1} of  $L_1(x,t)$ we have for $x$ in the integration interval and $t>0$ large
\[|L_1(x,t)|\ll 1+\frac{x^2}{t}+\frac{x^4}{t^2}.\]
Therefore,
\[I_1\ll t^{15/4}\int_{\alpha\log t}^\infty  e^{-\pi x/7}\,dx+
 t^{11/4}\int_{\alpha\log t}^\infty x^2 e^{-\pi x/7}\,dx+
t^{7/4}\int_{\alpha\log t}^\infty x^4 e^{-\pi x/7}\,dx.\]
Applying Gabcke's Lemma \ref{Gabcke}  
\[I_1\ll t^{-5/4}\log t+ t^{-9/4}(\log t)^3+t^{-13/4}(\log t)^5\ll t^{-5/4}\log t.\]
The second integral is 
\[
I_2:=e^{i\theta(t)}\rho_0(t)\int_{\alpha\log t}^\infty
L_1(x,t)\Bigl(\frac{t}{2\pi}\Bigr)^{ix/2}
\frac{\zeta(4+it+ix)}{7\cosh(\pi x/7)}\,dx\ll t^{15/4}
\int_{\alpha\log t}^\infty
\frac{|L_1(x,t)|}{\cosh(\pi x/7)}\,dx,
\]
which is bounded in the same way as $I_1$.
\end{proof}

\section{Second approximation to \texorpdfstring{$Z(t)$}{Z(t)}.}

\begin{proposition}
For all nonnegative integers $n$, we have
\begin{equation}\label{dercosh}
\frac{d^n}{dy^n} \frac{1}{\cosh y}=(-1)^n\frac{2e^{-y}B_n(-e^{-2y})}{(e^{-2y}+1)^{n+1}}.
\end{equation}
where $B_n(x)$ are polynomials of degree $n$.
\end{proposition}
\begin{proof}
For $n=0$ we have
\[\frac{1}{\cosh y}=\frac{2}{e^y+e^{-y}}=\frac{2e^{-y}}{1+e^{-2y}},\]
so that $B_0(x)=1$. 
By induction, we easily show that $B_n(x)$ satisfies the recurrence
\[B_n(x)=2x(1-x)B'_{n-1}(x)+(1+(2n-1)x)B_{n-1}(x).\]
From this equation it can be  shown that they are polynomials of degree $n$.
We will only need to know $B_0$, $B_1$ and $B_2$. The first few are
\[
\begin{array}{ll}
B_0(x)=1, &   B_3(x)=x^3+23 x^2+23x+1,    \\
B_1(x)=x+1,  & B_4(x)=x^4+76 x^3+230 x^2+76 x+1,    \\
B_2(x)=x^2+6 x+1,  & B_5(x)=x^5+237 x^4+1682 x^3+1682 x^2+237x+1.
\end{array}\qedhere
\]
\end{proof}
The coefficient of these polynomials appears in Oeis \oeis{A060187}. 
These polynomials are known as \emph{type B Eulerian polynomials} \cite{Pet}. They have 
the following  generating functions \cite{MR3408615}*{eq.~(13.3), (13.6)}
\[\sum_{j\ge0}x^j(2j+1)^n=\frac{B_n(x)}{(1-x)^{n+1}},\qquad \sum_{n=0}^\infty B_n(x)\frac{t^n}{n!}=\frac{(1-x)e^{(1-x)t}}{1-xe^{2(1-x)t}}.\]

We will use the following integral
\begin{proposition}\label{P:integral}
For $n$ a nonnegative integer, $\alpha$ and $t$ real numbers we have
\begin{equation}\label{E-171101-1}
\int_{-\infty}^\infty \frac{e^{i\alpha x} x^n}{7\cosh\frac{\pi x }{7}}\,dx=
\frac{2 e^{-y}B_n(-e^{-2y})}{(1+e^{-2y})^{n+1}}\Bigl(\frac{7i}{2}\Bigr)^{n},\qquad 
\text{with } y=\frac{7\alpha}{2}.
\end{equation}
\end{proposition}
\begin{proof}
It is known (\cite{SS}*{p.~81}) that for $y\in\R$
\[\int_{-\infty}^\infty \frac{e^{2\pi i y x}}{\cosh \pi x}\,dx=\frac{1}{\cosh\pi y}, 
\quad \text{so that}
\quad \int_{-\infty}^\infty \frac{e^{2 i y x}}{\cosh \pi x}\,dx=\frac{1}{\cosh  y}.\]
Differentiating, we get 
\[(2i )^n\int_{-\infty}^\infty \frac{x^n e^{2 i y x}}{\cosh \pi x}\,dx=
\frac{d^n}{dy^n} \frac{1}{\cosh y}=(-1)^n\frac{2e^{-y}B_n(-e^{-2y})}{(e^{-2y}+1)^{n+1}}.\]
Therefore, our integral can be calculated, first changing variables $x\mapsto 7x$ and then  taking $y=\frac{7\alpha}{2}$
\[\int_{-\infty}^\infty \frac{e^{i\alpha x} x^n}{7\cosh\frac{\pi x }{7}}\,dx=
7^n\int_{-\infty}^\infty \frac{e^{i7\alpha  x} x^n}{\cosh \pi x}\,dx=
\Bigl(-\frac{7}{2i}\Bigr)^n\frac{2e^{-y}B_n(-e^{-2y})}{(e^{-2y}+1)^{n+1}}.\qedhere\]
\end{proof}

We now define the auxiliary functions $H_r(t)$.
\begin{definition}
For $t$ real and $r$ a nonnegative integer we define
\[H_r(t)=\int_{-\infty}^\infty e^{i\beta x}
\frac{x^r\zeta(4+it+ix)}{7\cosh(\pi x/7)}\,dx,\qquad \beta=\frac12\log\frac{t}{2\pi},\]
\end{definition}
Then the definition of $F_4(t)=G(t)$ implies that 
\begin{multline}\label{E:Gexpr}
G(t)=e^{i\theta(t)}\rho_0(t)\cdot\\
\cdot\Bigl(H_0(t)+\frac{15H_1(t)}{4t}
+\frac{i H_2(t)}{4t}+\frac{165 H_2(t)}{32t^2}+\frac{241i H_1(t)}{24t^2}+\frac{41 i H_3(t)}{48t^2}-\frac{H_4(t)}{32t^2}\Bigr).
\end{multline}

\begin{proposition}\label{P:boundH}
When $t\to+\infty$ and a fixed nonnegative integer $r$ we have
$H_r(t)=\Orden(t^{-3/2})$.
\end{proposition}

\begin{proof}
Since the Dirichlet  series for $\zeta(4+it)$ converges uniformly in $\R$, we 
may interchange sum and integral, so that applying Proposition \ref{P:integral}
we obtain
\[
H_r(t)=\sum_{n=1}^\infty \frac{1}{n^{4+it}}
\int_{-\infty}^\infty e^{i(\beta-\log n) x}
\frac{x^r}{7\cosh\frac{\pi x}{7}}\,dx=\Bigl(\frac{7i}{2}\Bigr)^r
\sum_{n=1}^\infty \frac{1}{n^{4+it}}\frac{2e^{-y_n}B_r(-e^{-2y_n})}{(1+e^{-2y_n})^{r+1}},
\]
where
\[y_n=\frac{7}{2}(\beta-\log n)=\frac74\log\frac{t}{2\pi n^2}, \quad 
e^{-y_n}=\Bigl(\frac{2\pi n^2}{t}\Bigr)^{7/4}.\]

Since $B_r(x)$ is a polynomial of degree $r$. There is a constant $C_r$ such that 
\[B_r(-x)\le C_r(1+x)^r, \qquad x\ge0.\]
Then 
\begin{align*}
|H_r(t)|&=\Bigl|\Bigl(\frac{7}{2i}\Bigr)^r
\sum_{n=1}^\infty \frac{1}{n^{4+it}}\frac{2e^{-y_n}B_r(-e^{-2y_n})}{(1+e^{-2y_n})^{r+1}}\Bigr|
\le C_r\Bigl(\frac{7}{2}\Bigr)^r 
\sum_{n=1}^\infty \frac{1}{n^{4}}\frac{2e^{-y_n}}{(1+e^{-2y_n})}\\
&\le C_r\Bigl(\frac{7}{2}\Bigr)^r 
\sum_{n=1}^\infty \frac{1}{n^{4}}\frac{2}{\left(\frac{2\pi n^2}{t}\right)^{7/4}
+\left(\frac{2\pi n^2}{t}\right)^{-7/4}}\\&=
 C_r\Bigl(\frac{7}{2}\Bigr)^r 
\sum_{n=1}^\infty \frac{2}{\left(\frac{2\pi}{t}\right)^{7/4} n^{15/2}
+\left(\frac{2\pi}{t}\right)^{-7/4}n^{1/2}}\\
&\le C_r\Bigl(\frac{7}{2}\Bigr)^r
\int_{0}^\infty \frac{2\,dx}{\left(\frac{2\pi}{t}\right)^{7/4} x^{15/2}
+\left(\frac{2\pi}{t}\right)^{-7/4}x^{1/2}}.
\end{align*}
Change variables $x=(\frac{t}{2\pi})^{1/2}y$ and we obtain 
\[|H_r(t)|\le C_r\Bigl(\frac{7}{2}\Bigr)^r\Bigl(\frac{2\pi}{t}\Bigr)^{3/2}\int_0^\infty \frac{2dx}{y^{15/2}+y^{1/2}}\ll t^{-3/2}.\qedhere\]
\end{proof}

\begin{theorem}
For $t>0$ let 
\begin{equation}\label{D:H}
H(t)=H_0(t)=\int_{-\infty}^\infty\Bigl(\frac{t}{2\pi}\Bigr)^{ix/2}\frac{\zeta(4+it+ix)}{7\cosh(\pi x/7)}\,dx.\end{equation}
Then we have
\begin{equation}
Z(t)=\Bigl(\frac{t}{2\pi}\Bigr)^{\frac74}\Re\{e^{i\theta(t)}H(t)\}+\Orden(t^{-3/4}).
\end{equation}
\end{theorem}
\begin{proof} By \eqref{E:Gexpr} we may write \eqref{E:1staprox}, as 
\[Z(t)=\frac{\Re\{e^{i\theta(t)}\rho_0(t) H_0(t)\}}{((\frac14+t^2)(\frac{25}{4}+t^2))^{1/2}}+\sum_{j=1}^6 T_j+\Orden(t^{-5/4}\log^6t),\]
where each $T_j$ is the real part of is of the form  \[T_j =C_je^{i\theta(t)}\rho_0(t)((\tfrac14+t^2)(\tfrac{25}{4}+t^2))^{-1/2}H_r(t)t^{-k}\] where $C_j$ is a constant for each $j$, $r\in\{1,2,3\}$ and $k\in\{1,2\}$. Hence by Proposition \ref{P:boundH}
\[|T_j|\ll t^{2+\frac74}t^{-2}t^{-3/2} t^{-k}=t^{\frac14-k}\ll t^{-3/4}.\]
Since
\begin{align*}
\frac{\rho_0(t)}{((\frac14+t^2)(\frac{25}{4}+t^2))^{1/2}}&=(2\pi)^{-7/4}t^{2+\frac74}\Bigl(1+\frac{19}{t^2}+\Orden(t^{-4})\Bigr)\Bigl(\frac{1}{t^2}-\frac{13}{4t^4}+\Orden(t^{-6})\Bigr)\\
&=\Bigl(\frac{t}{2\pi}\Bigr)^{\frac74}+\Orden(t^{-\frac14}),
\end{align*}
we have,
\[Z(t)=\Bigl(\frac{t}{2\pi}\Bigr)^{\frac74}\Re\{e^{i\theta(t)}H(t)\}+\Orden(t^{-\frac14}H_0(t))+\Orden(t^{-3/4}).\]
Since $\Orden(t^{-\frac14}H_0(t))=\Orden(t^{-7/4})$ the result follows.
\end{proof}
To give an idea of the value of this approximation, I give a table of values. 

\begin{table}[H]
	\centering
	\begin{tabular}[t]{|l|l|l||l|l|l|}
		\hline
		& & & & & \\[-2ex]		
		$t$ & $Z(t)$ & $\Re\{(t/2\pi)^{7/4}e^{i\theta(t)}H(t)\}$  & $t$ & $Z(t)$ & $\Re\{(t/2\pi)^{7/4}e^{i\theta(t)}H(t)\}$ \\[0.4ex]
		\hline
		& & & & & \\[-2ex]	
		10 & $-1.5491945$ & $-0.9983260$ & $10^5$ & $\phantom{-}5.8795925$  & $\phantom{-}5.8790158$\\[0.4ex]
		\hline
		& & & & & \\[-2ex]
		$10^2$ & $\phantom{-}2.6926971$ & $\phantom{-}2.6269297$ & $10^6$ & $-2.8061339$ & $-2.8061012$\\[0.4ex]
		\hline
		& & & & & \\[-2ex]
		$10^3$ & $\phantom{-}0.9977946$ & $\phantom{-}0.9849027$ & $10^7$ & $\mskip7mu 14.3525504$ & $\mskip7mu 14.3525613$\\[0.4ex]
		\hline
		& & & & & \\[-2ex]
		$10^4$ & $-0.3413947$ & $-0.3452059$ & $10^8$ & $\phantom{-}3.6454079$ & $\phantom{-}3.6454066$\\[0.4ex]
		\hline
	\end{tabular}
\end{table}

\begin{proposition}
For $t>0$ we have 
\[H(t)=\sum_{n=1}^\infty \frac{1}{n^{4+it}}\frac{2}{(\frac{t}{2\pi n^2})^{7/4}+(\frac{t}{2\pi n^2})^{-7/4}}.\]
\end{proposition}
\begin{proof}
In the integral \eqref{D:H} defining $H(t)$ we expand $\zeta(4+it+ix)$ in series and integrate term by term. Then we apply \eqref{E-171101-1}. 
\begin{align*}
H(t)&=\int_{-\infty}^\infty\Bigl(\frac{t}{2\pi}\Bigr)^{ix/2}\frac{\zeta(4+it+ix)}{7\cosh(\pi x/7)}\,dx=\sum_{n=1}^\infty \frac{1}{n^{4+it}}\int_{-\infty}^\infty\Bigl(\frac{t}{2\pi n^2}\Bigr)^{ix/2}\frac{1}{7\cosh(\pi x/7)}\,dx\\
&=\sum_{n=1}^\infty \frac{1}{n^{4+it}}\frac{2}{(\frac{t}{2\pi n^2})^{7/4}+(\frac{t}{2\pi n^2})^{-7/4}}.\qedhere
\end{align*}
\end{proof}

\begin{remark}
The above expression can be written as 
\[H(t)=\Bigl(\frac{t}{2\pi}\Bigr)^{-7/4}\sum_{n=1}^\infty \frac{1}{n^{\frac12+it}}\frac{2}{1+(\frac{t}{2\pi n^2})^{-7/2}}.\]
\end{remark}

\begin{remark}
It appears that the phase of $H(t)$ decreases moderately with increasing $t$. For example, in the figure, we see that $\arg H(t)$ for  $t\in (10000,10020)$ start in the fourth quadrant and purely imaginary (or almost there), pass to the third quadrant, first quadrant and then 

\begin{figure}[H]
\begin{center}
\includegraphics[width=\hsize]{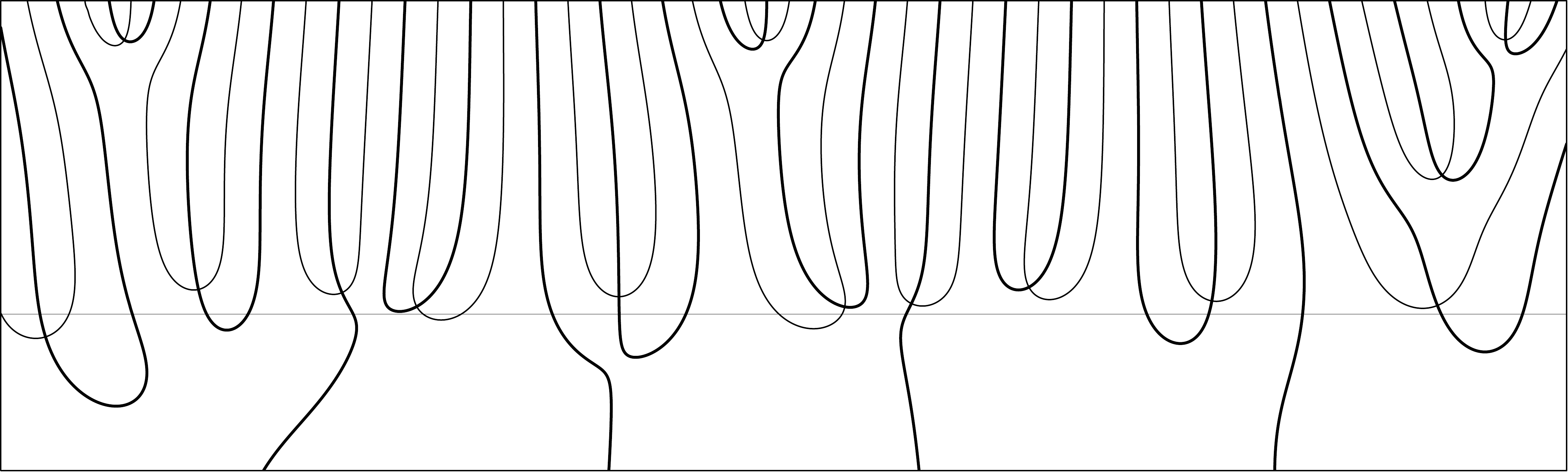}
\caption{x-ray of $H(z)$ for $z\in(10000,10020)\times(-2,4)$}
\label{Happrox}
\end{center}
\end{figure}
\noindent move between 1st and 4th quadrant (with two slight changes into 2nd and 3rd) ending in the fourth quadrant. Therefore, we pass from $-\pi/2$ to some point in the interval $(-2\pi, -2\pi-\pi/2)$. So the variation of the argument is at most $-2\pi$.  In this rectangle there is only one zero of $H(z)$ with $\Im z<0$. 
\end{remark}

\begin{remark}
The computation of the continuous argument of $H(t)$ for $0<t<15000$ is consistent with 
\[0\le -\arg H(t)\approx c\Bigl(\frac{t}{2}\log\frac{t}{2\pi}-\frac{t}{2}\Bigr),\qquad c=\frac{17}{63}\approx0.27.\]
If this continues to be true for larger values of $t$, it will imply 73\% of zeros of zeta on the critical line.
\end{remark}

\section*{Appendix. Counting zeros by means of the phase.}\label{appendix}
We show here how finding real analytic functions $F(t)$ such that $Z(t)=\Re F(t)$ is relevant to counting zeros on the critical line.

Let $f\colon[a,b]\to\C$. We say that it is \emph{real analytic} if $f$ has an analytic 
extension to a neighborhood of $[a,b]$ in $\C$.  Notice that $f$ real analytic implies
that $g(t)=\overline{f(t)}=\overline{f(\overline{t})}$ is real analytic also, so that 
\[\Re f(t)=\frac{f(t)+g(t)}{2},\qquad \Im f(t)=\frac{f(t)-g(t)}{2i},\]
are also real analytic functions in $[a,b]$.

First we extend a Proposition from \cite{AL}*{p.~218} 
\begin{proposition} \label{P:1}
If $f\colon[a,b]\to\C$ is real analytic not identically $0$, then there are two real 
analytic functions $U\colon[a,b]\to\R$ and $\varphi\colon[a,b]\to\R$ such that
\[f(t)=U(t)e^{i\varphi(t)}.\]
Given two such representations $f=U_1e^{i\varphi_1}$
and $f=U_2e^{i\varphi_2}$, we have either $U_1=U_2$ and $\varphi_1-\varphi_2=2k\pi$
or $U_1=-U_2$ and $\varphi_1-\varphi_2=(2k+1)\pi$ for some integer $k$.
\end{proposition}

\begin{proof}
If $f(t)\ne0$ there is a simply connected neighborhood $\Omega$ of $[a,b]$ where
$f$ extends as an analytic function that is never $0$ in $\Omega$. Then we may 
define an analytic function  $g\colon\Omega\to\C$ such that $e^{g(t)}=f(t)$
Then we may take $U(t)=e^{\Re g(t)}$ and $\alpha(t)=\Im g(t)$. 

If $f$ have zeros in $[a,b]$, they are a finite number (counting multiplicities).
Let $\alpha_j$, with $1\le j\le n$ these zeros. Let 
$F(t)=f(t)\prod_{j=1}^n (t-\alpha_j)^{-1}$. Then $F$ is real analytic and does not 
vanish, so it can be written as  $F(t)=V(t)e^{i\alpha(t)}$. Putting
$U(t)=V(t)\prod_j (t-\alpha_j)$ we obtain our desired representation.

Given two representations $f=U_1 e^{i\varphi_1}=U_2e^{i\varphi_2}$, $U_2$ is not identically $0$, so  we will 
have $U_1/U_2=e^{i(\varphi_2-\varphi_1)}$. Therefore
$U_1/U_2$ is real analytic, take only real values, and have modulus $=1$. So 
$U_1/U_2=1$ or it is equal to $-1$.  Then the phases have to take the given 
values.
\end{proof}

Given such a representation $f=Ue^{i\varphi}$ we say that $\varphi$ is a phase of $f$.
When $f(t)$ do not vanish in $[a,b]$ we can take $U(t)=|f(t)|$ and then 
for any $a\le t\le b$, the value of  $\varphi(t)$ coincides with some argument of 
$f(t)$.

Given a real analytic function $f\colon[a,b]\to \C$, not identically $0$, we define 
$N_f[a,b]$ as the number of zeros of $f$ in the closed segment $[a,b]$
counting multiplicities. 

\begin{proposition}\label{P:171207-2}
Let $f\colon[a,b]\to\R$ and $g\colon[a,b]\to \C$ be real analytic functions
with $f(t)=\Re g(t)$ not identically $0$. Let $\varphi(t)$ a phase of $g(t)$. Then 
\[\frac{|\varphi(b)-\varphi(a)|}{\pi}<  N_f(a,b)+1.\]
\end{proposition}

\begin{proof}
Substituting $g$ by $\overline{g}$ that is also real analytic and 
has the same real part, we may assume that $\varphi(a)\le \varphi(b)$.

We have  something to prove only in the case  $\varphi(b)-\varphi(a)> \pi$.
Consider the odd integer multiples of $\pi/2$ contained  between $\varphi(a)$ 
and $\varphi(b)$. Since $\varphi(b)-\varphi(a)> \pi$ there is at least one. 
Therefore, there exist $n\le m$ integers such that 
\[
(2n-1)\frac{\pi}{2}<\varphi(a)\le (2n+1)\frac{\pi}{2}<(2n+3)\frac{\pi}{2}<
\cdots<
(2m+1)\frac{\pi}{2}\le\varphi(b)<(2m+3)\frac{\pi}{2}.
\]
Since $\varphi(t)$ is continuous, it takes all values in the 
interval $[\varphi(a),\varphi(b)]$. Therefore, there are points $t_k\in[a,b]$ with 
$\varphi(t_k)=(2k+1)\frac{\pi}{2}$ for  $n\le k\le m$. At these points, we have
$f(t_k)=\Re g(t_k)=0$. That is $N_f(a,b)\ge m-n+1$.  The 
interval $((2n-1)\pi/2,(2m+3)\pi/2)\supset[\varphi(a),\varphi(b)]$, comparing 
its lengths, yields 
\[\pi(m-n)+2\pi> \varphi(b)-\varphi(a).\]
The result follows.
\end{proof}

\begin{proposition}
Let $f$ and $g\colon[a,b]\to\C$ real analytic functions such that 
\begin{equation}\label{E:171207-1}
|f(t)-g(t)|<|f(t)|, \qquad \text{for all } a\le t\le b.
\end{equation}
Then $f$ and $g$ do not vanish on $[a,b]$, so that we can write 
$f(t)=|f(t)|e^{i\alpha(t)}$ and $g(t)=|g(t)|e^{i\beta(t)}$ with $\alpha$ and $\beta$
real analytic.  Then we can choose $\beta(t)$ so that 
\[ |\alpha(t)-\beta(t)|<\pi,\qquad a\le t\le b.\]
In particular 
\[\Bigl|\frac{\alpha(b)-\alpha(a)}{\pi}- \frac{\beta(b)-\beta(a)}{\pi}\Bigr|<2.\]
\end{proposition}

\begin{proof}
If $f(t)=0$, equation \eqref{E:171207-1} implies $|g(t)|<0$ which is a contradiction.
If $g(t)=0$, then we arrive at the contradiction $|f(t)|<|f(t)|$.  Since the functions
do not vanish, we may apply Proposition \ref{P:171207-2} to get
$f(t)=|f(t)|e^{i\alpha(t)}$ and $g(t)=|g(t)|e^{i\beta(t)}$.  The function $\beta(t)$ is 
not unique, in fact $\beta(t)+2k\pi$ for any integer $k$ can be used instead of 
$\beta(t)$. 

We have $|f(0)-g(0)|<|f(0)|$. Therefore $g(0)$ is in the open disc with center $f(0)$ 
and radius the distance from $f(0)$ to $0$.  $\alpha(0)$ is one of the determination 
of $\arg f(0)$, there is a value of $\arg g(0)$ such that $|\alpha(0)-\arg g(0)|<\pi$, since
the two points are in some half plane.  We may choose a $k$ so that 
$|\beta(0)+2k\pi-\alpha(0)|<\pi$, without loss of generality we may assume that 
$|\beta(0)-\alpha(0)|<\pi$. Since these functions are continuous, this will be true 
for some interval $|\beta(t)-\alpha(t)|<\pi$ if $a\le t<a+\varepsilon$. 

By contradiction we shall prove that $|\beta(t)-\alpha(t)|<\pi$ for all $t\in [a,b]$. 
Assume that there is some $|\beta(t)-\alpha(t)|\ge\pi$. The set of points that satisfies
this inequality is a closed set of $[a,b]$ contained in $[a+\varepsilon, b]$.
Let $t_0$ be the infimum of these closed sets. By continuity we will have 
$|\beta(t_0)-\alpha(t_0)|=\pi$. Therefore $f(t_0)=|f(t_0)|e^{i\alpha(t_0)}$ and
\[g(t_0)=|g(t_0)|e^{i\beta(t_0)}=
|g(t_0)|e^{i\alpha(t_0)\pm \pi i}=-|g(t_0)|e^{i\alpha(t_0)}.\]
And we have a contradiction with \eqref{E:171207-1}
because
\[|f(t_0)|>|f(t_0)-g(t_0)|=\bigl||f(t_0)|e^{i\alpha(t_0)}+|g(t_0)|e^{i\alpha(t_0)}\bigr|=
|f(t_0)|+|g(t_0)|.\]
The rest is an easy consequence of this.
\end{proof}

\end{document}